\documentclass[11pt]{amsart}

\usepackage{latexsym,amssymb,amsmath,graphics}
\usepackage[dvips]{graphicx}

\def\ffi{\varphi}
\def\eps{\varepsilon}
\def\dst{\displaystyle}

\def\supp{{\mathrm{supp}\,}}

\def\N{{\mathbb{N}}}

\def\Q{{\mathbb{Q}}}
\def\R{{\mathbb{R}}}

\def\Z{{\mathbb{Z}}}

\newcommand{\norm}[1]{{\left\|{#1}\right\|}}
\newcommand{\ent}[1]{{\left[{#1}\right]}}
\newcommand{\abs}[1]{{\left|{#1}\right|}}
\newcommand{\scal}[1]{{\left\langle{#1}\right\rangle}}

\newenvironment{notation}[1][]{\vskip1pt\noindent\rm\textit{Notation}\,:\ }{\rm\vskip1pt}
\newenvironment{remark}[1][]{\vskip1pt\noindent\rm\textit{Remark #1}\,:\ }{\rm\vskip1pt}

\newtheorem{problem}{Problem}

\newtheorem{lemma}{Lemma}[section]
\newtheorem{proposition}[lemma]{Proposition}
\newtheorem{theorem}[lemma]{Theorem}

\begin{document}

\title{Reconstruction of functions from their triple correlations}
\author{Philippe Jaming \& Mihail N. Kolountzakis}
\address{P.J.\,: Universit\'e d'Orl\'eans\\
Facult\'e des Sciences\\ 
D\'epartement de Ma\-th\'e\-ma\-ti\-ques\\BP 6759\\ F 45067 Orl\'eans Cedex 2\\
France}
\email{jaming@labomath.univ-orleans.fr}

\address{M.K.\,: Department of Mathematics\\ University of Crete\\ Knossos Ave.\\
714 09 Iraklio\\ Greece}
\email{kolount@member.ams.org}

\begin{abstract}
Suppose that $A$ is a subset of an abelian group $G$.
To know the {\em 3-deck} of $A$ is to know the number of
occurrences in $A$ of translates of each possible
multiset $\{0,a,b\}$.
The concept of the 3-deck of a set is naturally extended to
$L^1$ functions on $G$.
In this paper we study when the 3-deck of a function 
determines the function up to translations.
The method is to look at the Fourier Transform of the function.
Our emphasis is on the real line and the cyclic groups.
\end{abstract}

\date{June 20, 2003}

\thanks{Research partially financed by : {\it European Commission}
Harmonic Analysis and Related Problems 2002-2006 IHP Network
(Contract Number: HPRN-CT-2001-00273 - HARP).\\
The authors wish to particularly thank
Professor C.C. Moore for providing several useful references.}

\maketitle

\section{Introduction}

The aim of this article is to address a problem that arises independently in
combinatorics and in diffraction theory.

In combinatorics, a common problem is to reconstruct an object (up to
isomorphism) from the collection of isomorphism classes of its sub-objects.
The most famous problems of that nature are the {\it Reconstruction conjecture}
and the {\it Edge reconstruction conjecture} and we refer to Bondy \cite{Bo}
and Bondy-Hemminger \cite{BH} for surveys on the subject.
In that direction, Radcliffe and Scott investigated in \cite{RS1,RS2} the problem of
reconstructing subsets of $\R$, $\Z^n$, $\Q^n$ or $\Z/n\Z$ from the collection
of isomorphism classes of their subsets of prescribed size. (Here two sets
are isomorphic if they are translates of each other).

To be more precise, fix an integer $k\geq3$ and
let $G$ be one of the groups $\R$, $\Z^n$, $\Q^n$ or
$\Z/n\Z$.
A subset $E$ of $G$
is {\it locally finite} if $E\cap(E-x)$ is finite for every $x\in G$. The 
{\it $k$-deck} of $E$ is then the function on $G^{k-1}$ defined by
$$
N_E(x_1,\ldots,x_{k-1})=\abs{E\cap(E-x_1)\cap\cdots\cap(E-x_{k-1})}
$$
where $\abs{.}$ stands for cardinal.

The problem addressed by Radcliffe and Scott is the following:
\begin{problem}
\label{prb:main}
Assume that $E$ and $F$ are two locally finite sets
in an abelian group $G$ that
have the same $k$-deck, for some integer $k$.
Are $E$ and $F$ translates of each other?
\end{problem}
It turns out that, in the particular case of $G=\Z/n\Z$,
this problem also arises in the mathematical theory of diffraction.

It is well known that the answer is negative for $k=2$, but
for $k\geq 4$, Grunbaum and Moore showed that the answer is positive.
This result can be further improved in some cases.
For example
it is shown in \cite{RS1} that if $n$ is a prime, then the answer is still positive
for $k=3$.

One of our aims here is to summarize the knowledge on the problem in the $\Z/n\Z$ case,
to improve some of the proofs and to give the last possible positive results in the cases
$n=p^a$ and $n=pq$ with $p,q$ primes. This bridges the results in \cite{RS1} and \cite{GM},
and essentially closes the problem.
To do so, we make a more explicit use of the Fourier transform than in \cite{RS1}.

Radcliffe and Scott in \cite{RS2} further prove that the answer to Problem \ref{prb:main} is
positive if the group is one of $\R$, $\Z^n$, $\Q^n$. 
Moreover they ask for measure-theoretic  counterparts to their results in $\R$.
This is the main question that we will address here.

To be more precise, we will consider the Lebesgue measure on $\R$
and write again $\abs{E}$ for the measure of a measurable subset $E$ of $\R$.
The {\it $k$-deck} of $E$ is still defined as the function on $\R^{k-1}$ given by
$$
N_E(x_1,\ldots,x_{k-1})=\abs{E\cap(E-x_1)\cap\cdots\cap(E-x_{k-1})}.
$$
We will also consider naturally the $k$-deck of a nonnegative
function $f\in L^1(\R)$ defined by
$$
N_f(x_1,\ldots,x_{k-1})=\int f(t)f(t+x_1)\cdots f(t+x_{k-1})dt.
$$
If $\chi_E$ is the characteristic function of a set $E$ of finite measure, then
$N_E=N_{\chi_E}$. 
The problem we address is the following:

\begin{problem}
\label{prb:measure-theoretic}
Let $k\ge 3$ be an integer and let $N_f$ denote the $k$-deck of the function $f$.

(a) Suppose $0\le f,g\in L^1(\R)$ such that $N_f=N_g$ almost everywhere.
Does there exist an $a\in\R$ such that $g(x)=f(x-a)$ for almost every
$x\in\R$\,?

(b) Let $E,F$ be two sets of finite measure such that $N_E=N_F$ almost everywhere.
Does there exist an $a\in\R$ such that $F=E-a$ up to a set of
zero measure.

If this happens we say that $f$ (resp.\ $E$) is \emph{determined up to
translation by its $k$-deck}.
\end{problem}

Again, this problem occurred before in the context of texture analysis as well as
in crystallography where the
$k$-deck is known as an {\em higher order autocorrelation function} ({\it see}
\cite{AK,CW,Rot}). In these fields, a common problem is to be
able to reconstruct a function from its $2$-deck or autocorrelation function.
This is easily seen to be equivalent
to that of reconstructing a function knowing only the modulus of its Fourier
Transform and is therefore called the {\em phase retrieval problem}.
We refer to the surveys \cite{BM,KST,Mi}
the book \cite{Hu} or the introduction of \cite{Ja} for further references.
Note that the $3$-deck has been proposed to overcome the non-uniqueness in
the reconstruction ({\it see} \cite{Mi}).

We will only give a partial answer to this question and show that in many cases,
the $3$-deck will be enough for the answer to be yes. In particular, we managed
to do so for characteristic functions of compact sets, a fact independently
proved by Rautenbach and Triesch \cite{RT}.
These positive results are described in Section \ref{sec:positive-results-on-R}.

In the opposite direction, for $k$ fixed, building on a construction in \cite{Ja}
(and on an earlier construction in \cite{CW}),
we construct in Section \ref{sec:negative-results-on-R} an uncountable family
of functions $(f_i)$ that all have same $3$-deck but are not translates of one another.
However, these functions are not characteristic functions of sets,
and we have not been able to settle completely the second problem.
Nevertheless, in Section \ref{sec:stability}, we have managed to prove that
if $f=\chi_E$ and $g\in L^1(\R)\cap L^\infty(R)$ is non-negative, then 
$N_g=N_f$ implies that $g$ is also a characteristic function of some set $F$,
thus reducing Problem \ref{prb:measure-theoretic} (a) to Problem \ref{prb:measure-theoretic} (b).
However we do not know if $F=E-a$ for some $a\in\R$.

In Section \ref{sec:restricted} we show that if $E$ belongs to a rather wide class
of measurable sets of finite measure then it is determined up to translation from its 3-deck, among all measurable sets.
This class consists of all sets which are 
unions of intervals with endpoints that do not have an accumulation point
and such that the complementary intervals have length bounded below by a positive constant.

To conclude this introduction, let us point out that our investigation on the
$3$-deck problem for sets establishes a link between phase retrieval problems
and uncertainty principles. Namely, the obstruction to uniqueness
in the $3$-deck problem is the possibility to have large gaps in the spectrum
({\it i.e.} support of the Fourier Transform) of a function.
Although this is easy to
achieve for general functions, it is not known whether characteristic functions
may have many gaps in the support of their Fourier Transform.
The best results in that direction, due to Kargaev and Volberg \cite{Ka,KV},
state only that the Fourier Transform of a characteristic function may be
zero on an interval.
However this is not enough for proving a counterexample (to determination from 3-deck)
as the results of Section \ref{sec:restricted} are applicable to sets of the type
used by Kargaev and Volberg, where we prove that the specific sets used in \cite{Ka,KV}
are determined by there $3$-deck up to translation.
\vskip6pt
This article is organized as follows.
In the Section \ref{sec:results-on-R}, we will focus on results on the real line, starting with the reformulation
of the problem in terms of the Fourier Transform, continuing with positive and negative results
and concluding with the stability result mentioned above.
In Section \ref{sec:cyclic}
we focus on the $\Z/n\Z$ case, starting with a survey of known results with the aim
of bridging the knowledge in both communities that are interested in the problem.
We are then focusing on our own results.

\section{The results on $\R$}
\label{sec:results-on-R}

\subsection{Reformulation in terms of Fourier Transforms}
\label{reformulation}

In this section we will reformulate the problem in terms of the Fourier
Transform. For the the sake of simplicity of notation, we will only focus on the case
of the real line, although the results generalize without difficulty to
other locally compact abelian groups.

Also the results here have been proved in a similar way in \cite{CW} and again in \cite{RT}. As we
need the notations, we prefer reproducing the proofs for sake of simplicity
and completeness.

To fix notation, we define the Fourier Transform of
$\ffi\in L^1(\R^d)$  by
$$
\hat\ffi(\xi)=\int_{\R^d}\ffi(t)e^{-2\pi i\scal{t,\xi}}dt.
$$
It is easy to see that if $f\in L^1(\R)$, then $N_f\in L^1(\R^{k-1})$
with $\norm{N_f}_1\leq\norm{f}_1^k$, with equality when $f$ is non-negative. 
It is not difficult to see that actually $\norm{N_f}_r\leq\prod_{j=1}^k\norm{f}_{p_j}$
as soon as $1+\frac{k-1}{r}=\sum_{j=1}^k\frac{1}{p_j}$ and that $N_f$ is continuous
as soon as $f\in L^1(\R)\cap L^p(\R)$ for some $p\geq3$. Actually all
results for convolutions apply for $k$-decks with essentially the same proofs.

Further, the Fourier Transform of
$N_f$ is
$$
\widehat{N_f}(\xi_1,\ldots,\xi_{k-1})=\widehat{f}(\xi_1)\cdots\widehat{f}(\xi_{k-1})
\widehat{f}\bigl(-(\xi_1+\cdots+\xi_{k-1})\bigr).
$$
This implies that solving $N_g=N_f$ is equivalent to solving
\begin{equation}
\label{eq1}
\widehat{g}(\xi_1)\cdots\widehat{g}(\xi_{k-1})
\widehat{g}\bigl(-(\xi_1+\cdots+\xi_{k-1})\bigr)
=\widehat{f}(\xi_1)\cdots\widehat{f}(\xi_{k-1})
\widehat{f}\bigl(-(\xi_1+\cdots+\xi_{k-1})\bigr)
\end{equation}
for all $\xi_1,\ldots,\xi_{k-1}\in\R$.

As our primary interest is in characteristic functions, we will from now on
consider that $f$ and $g$ are {\it nonnegative $L^1$ functions} and that
$\widehat{f}(0):=\int_\R f(t)dt\not=0$. 
Another reason for considering only real functions is to avoid the extra
complication due to the fact that the functions $f$ and $\omega f$ have the
same $k$-deck whenever $\omega$ is a $k$-th root of unity.

\medskip

In the sequel {\it we will only consider Problem \ref{prb:measure-theoretic}
for non-negative functions}.

\medskip

Then, taking $\xi_1=\cdots=\xi_{k-1}=0$ in (\ref{eq1}) gives
$\widehat{g}(0)^k=\widehat{f}(0)^k$, and since $\widehat{g}(0)$
and $\widehat{f}(0)$ are nonnegative, they are equal.
Further, as $f$ and $g$ are real
$\widehat{f}(-\xi)=\overline{\widehat{f}(\xi)}$,

It follows that, if we take
$\xi_1=\xi$, $\xi_2=-\xi$ and $\xi_3=\cdots=\xi_{k-1}=0$ in (\ref{eq1}), we get 
$\abs{\widehat{g}(\xi)}=|\widehat{f}(\xi)|$.

Now, let us write
\begin{equation}
\label{eq:relation}
\widehat{g}(\xi)=e^{2\pi i \ffi(\xi)}\widehat{f}(\xi)
\end{equation}
where, by the continuity of $\widehat{f},\widehat{g}$, we may assume that $\ffi$ is
continuous on the support of $\widehat{f}$. Introducing this in equation
(\ref{eq1}), it reduces to the following
\begin{equation}
\label{eq2}
\ffi(\xi_1+\cdots+\xi_{k-1})=\ffi(\xi_1)+\cdots+\ffi(\xi_{k-1})\qquad(\bmod\ 1)
\end{equation}
for all $\xi_1,\ldots,\xi_{k-1}\in\supp\widehat{f}$
such that $\xi_1+\cdots+\xi_{k-1} \in \supp\widehat{f}$.

The solution of such an equation is then dependent on
the group on which one might consider it. In the case of the real line,
it is easy to show that the following holds
(see e.g.\ \cite{Ja}, Lemma 3):
\begin{lemma}
\label{lm:basic}
Suppose that $f, g \in L^1(\R)$ are nonnegative and have $N_f = N_g$.
It follows that $\widehat{f}, \widehat{g}$ are connected by \eqref{eq:relation}.
Write $\dst\mbox{supp}\,\widehat{f}=\bigcup_{j\in\Z}I_j$
to be the decomposition of $\mbox{supp}\,\widehat{f}$ into
disjoint open intervals numbered so that $0\in I_0$, $I_{-j}=-I_j$.
Then, there exists $\omega\in\R$ and a sequence $\theta_j$ that satisfies
$\theta_0=0$ and $\theta_{-j}=-\theta_j$ such that, if $x\in I_j$,
$$
\ffi(x)=\omega x+\theta_j.
$$
Moreover, if there are $\xi_1,\ldots,\xi_{k-1}$ with $\xi_l$ in some $I_{j_l}$ such that $\xi_1+\cdots\xi_{k-1}\in I_l$ then
\begin{equation}
\label{eq3}
\theta_{j_1}+\theta_{j_{k-1}}=\theta_l.
\end{equation}
In particular, if $\theta_j\not=\theta_l$ then $I_j$ and $I_l$ are distant by at
least $\frac{k-2}{2}\times\mathrm{diam}\,I_0$.
\label{fundlemma}
\end{lemma}

\begin{proof}[Sketch of proof] One first proves the result on $I_0\cap\Q$
where it is trivial, then extends it to all of $I_0$ by continuity.

The extension to $I_k$ is then obvious taking $x\in I_k,y\in I_0$. Finally
(\ref{eq3}) is a reformulation of (\ref{eq2}).
\end{proof}

\subsection{Positive results}
\label{sec:positive-results-on-R}

Recall that a measure $\mu$ has {\it divergent logarithmic integral} if one of 
the two following integrals
$$
\int_{-\infty}^0 \frac{-\log\abs{\mu}((-\infty, x])}{1 + x^2}~dx,
\ \ \ 
\int_0^\infty \frac{-\log\abs{\mu}([x, \infty))}{1 + x^2}~dx,
$$
is divergent.

We will now prove the following\,:

\begin{theorem}
Assume that $f\in L^1(\R)$ is real valued. 
In each of the following cases, the function $f$ is determined
up to translations by its $3$-deck\,:
\begin{enumerate}
\item\label{thpi} $f$ is of compact support;

\item\label{thpii} for some $M>0$ the integral $\int|f(x)|e^{M|x|}~dx$ is finite;

\item\label{thpiii} $\widehat{f}$ is analytic in a neighborhood of the real
line;

\item\label{thpiv} the measure $f\mbox{d}x$ has divergent logarithmic integral.
\end{enumerate}
\end{theorem}

\begin{proof} The case \ref{thpi} is a particular case of \ref{thpii} which in
turn is a particular case of \ref{thpiii}.
In this cases, $\{\widehat{f}=0\}$
is a discrete set.
In the last case, a theorem of Beurling (see \cite{Ko} p.\ 268) states that 
$\widehat{f}$ cannot be zero on a set of positive measure. In conclusion, this 
theorem is proved once we have proved the following lemma\,:

\begin{lemma} If $\widehat{f}$ does not vanish on a set of positive measure,
then $f$ is uniquely determined up to translations by its $3$-deck.
\label{lem:rk3}
\end{lemma}

Following the notation of Lemma \ref{fundlemma}, write 
$\dst\mbox{supp}\,\widehat{f}=\mbox{supp}\,\widehat{g}=\bigcup_{j\in\Z}I_j$ 
and let $\ffi(x)=\omega x+\theta_j$ be the function given by the lemma.
Write $I_0=(-a,a)$. Up to changing $g$ into $g(.-\omega)$, we may assume that
$\widehat{g}(\xi)=e^{2\pi i \theta_j}\widehat{f}(\xi)$ if $x\in I_j$.

Contrary to what we want to prove,
assume that $f$ is not almost everywhere equal to $g$ and
let $s = \inf\{x>0: \widehat{f}(x) \neq \widehat{g}(x)\}$.
It follows that $s$ is finite, otherwise $f\equiv g$.

From this we get that there is $j \in \N$ with $\theta_j\neq 0$
such that $I_j$ is of the form $I_j=]s,s+\lambda[$. 
By the last assertion of Lemma \ref{fundlemma}
we get that
$$
(0, s) \cap \{x>0:\ 0 \neq \widehat{f}(x) = \widehat{g}(x)\}
 \subseteq (0, s-a),
$$
which implies that $\widehat{f}\equiv 0$ on $(s-a, s)$,
in contradiction to our assumption on the support of $\widehat{f}$.
\end{proof} 

\begin{remark}
\begin{enumerate}
\item The proof of the lemma actually tells us that if
$\supp\widehat{f}$ has no gap of size $\mbox{diam}\,I_0$ then $f$ is determined
up to translations by its $3$-deck. This, as well as Point \ref{thpi} of the
theorem has been independently proved by Rautenbach and Triesch in \cite{RT}.

\item Further remarks of that order can be made.
For instance, one might notice that
$\mbox{diam}\,I_j$ has to stay bounded. Indeed, if this is not the case, we
may write $\xi\in I_j$ as $\xi=\xi_l+\xi-\xi_l$ with $\xi_l,\xi_l+\xi$ in some
$I_l$. Then $\theta_j=\ffi(\xi)=\ffi(\xi_l+\xi)+\ffi(-\xi_l)=0$.

\item As a corollary of the proof, we immediately get that if ${}^c\supp\widehat{f}$
is negligible, then the only solution of (\ref{eq3}) are $\ffi(\xi)=\omega\xi$ for
some $\omega\in\R$.
\end{enumerate}
\end{remark}

\subsection{Negative results}
\label{sec:negative-results-on-R}

\begin{theorem} For every $k\ge 3$, there exist two nonnegative and smooth functions
that have the same $k$-deck but are not translates of each other.
\end{theorem}

\begin{proof} We will only give full details of the proof in the case $k=3$.

Let $\psi(x)=\left(\frac{\sin \left(\frac{\pi}{2} x\right)}{\pi x}\right)^2$, then 
$\psi\in L^1(\R)$
and its Fourier Transform is
$$
\widehat{\psi}(\xi)=\begin{cases}
\left(\frac{1}{2}-\abs{\xi}\right)&\mbox{if
}\abs{\xi}\leq \frac{1}{2}\\ 0&\mbox{else}\\ \end{cases}.
$$ 

Consider $f(x)=\bigl(1+\cos(4\pi x)\bigr)\psi(x)$ so that
$\widehat{f}(\xi)=\widehat{\psi}(\xi+2)+\widehat{\psi}(\xi)+
\widehat{\psi}(\xi-2)$.
The support of $\widehat{f}$ is then $\ent{-\frac{1}{2},\frac{1}{2}}+\{-2,0,2\}$. 

Now define $\ffi$ to be such that $\ffi=-1$ on 
$\ent{-\frac{1}{2},\frac{1}{2}}+\{-2,2\}$, $\ffi=1$ on
$\ent{-\frac{1}{2},\frac{1}{2}}$ and $\ffi=0$ elsewhere. 
It is easy to see that $\ffi$ satisfies
$\ffi(x+y)=\ffi(x)\ffi(y)$ whenever $x,y,x+y\in\mbox{supp}\,\widehat{f}$.

Finally, let $g(x)=\bigl(1-\cos(4\pi x)\bigr)\psi(x)$ so that
$\widehat{g}(\xi)=\ffi(\xi)\widehat{f}(\xi)$.
Then $N_f=N_g$ but $f$ and $g$ are not translates of each other.

To adapt the proof to the case $k\geq4$, it is then enough to replace 
the $\cos(4\pi x)$ in the above argument by (say)
$\cos(2\frac{k+1}{2}\pi x)$.
\end{proof}

Using the same ideas as for the proof of Theorem 3 in \cite{Ja},
one may improve this to get the following result:

\begin{proposition} For every $k\ge 3$,
there exists a function $f\in L^1(\R)$ such that there are uncountably
many functions, not translates of each other,
that have the same $k$-deck as $f$.
\end{proposition}

\begin{proof}[Sketch of proof] The only thing to be done is to replace
the factor $1\pm\cos(k+1)\pi x$ by the Riesz product
$$
\prod_{n\in\Z}\left(1+a_n\eps_n\cos(k+1)^n\pi x\right)
$$
where $a_n$ decreases fast enough to ensure the convergence of the product and each
$\eps_n=\pm1$.
\end{proof}

On the other hand we cannot have this situation for all $k$:

\begin{lemma} If $f$ and $g$ are such that there $k$-decks are the same 
for every $k$, then $f$ and $g$ are translates of each other.
\end{lemma}

This is Theorem 2 in \cite{CW} (see also \cite{Rot} for some generalizations).
However, the proof bellow is a bit simpler.

\begin{proof}
Using the notation of Lemma
\ref{fundlemma}, let us write 
$\widehat{g}(\xi)=e^{2\pi i(\omega\xi+\theta_j)}\widehat{f}(\xi)$
if $\xi\in I_j$. Recall that the $\theta_j$'s satisfy, for every $k$,
\begin{equation}
\label{eqkdeckunique}
\theta_{j_1}+\cdots+\theta_{j_{k-1}}=\theta_l
\end{equation}
if there exists $\xi_1\in I_{j_1},\cdots,\xi_{k-1}\in I_{j_{k-1}}$
such that $\xi_1+\cdots\xi_{k-1}\in I_l$.
Write $I_l=(a_l,b_l)$.
Then, for $k\geq 2 \frac{a_l}{\mbox{diam}\,I_0}$, there exists 
$\xi_1,\cdots,\xi_{k-1}\in I_0$
such that $\xi_1+\cdots+\xi_{k-1}\in I_l$. Equation (\ref{eqkdeckunique})
then implies that $\theta_l=0$ so that $\widehat{g}(\xi)=e^{2\pi i\omega\xi}\widehat{f}(\xi)$
and $g$ is a translate of $f$. 
\end{proof}

\subsection{A stability result for characteristic functions}
\label{sec:stability}

In this section, we will prove that, under some mild restrictions,
if a function has the same $3$-deck as a set, then this function is a 
characteristic function of a set. More precisely, we will prove\,:

\begin{proposition} Let $E\subset\R$ be a set of finite measure, let $f=\chi_E$
and let $g\in L^1(\R)$ be a non-negative function such that $N_g=N_f$ almost everywhere.
Then there exists 
a set $F$ of finite measure such that $g=\chi_F$.
\end{proposition}

\begin{proof}
First note that, as $f=\chi_E$, $\norm{f}_1=\norm{f}_2^2$ and
$\norm{f}_3=\norm{f}_2^{2/3}$. Further, in this case $N_f$ is continuous.

Next, as we assumed that $g\geq0$, we get
$$
\norm{g}_1=\norm{N_g}_1^{1/3}=\norm{N_f}_1^{1/3}=\norm{f}_1.
$$
Further, as $\abs{\widehat{g}}=\abs{\widehat{f}}$ we get
$$
\norm{g}_2=\norm{\widehat{g}}_2=\norm{\widehat{f}}_2=\norm{f}_2
$$
with Parseval. 

Now, if $g\notin L^3(\R)$ then, by positivity of $g$, $N_g$ is not essentially bounded in a
neighborhood of $0$.
This would contradict $N_g=N_f$ almost everywhere.

Now that we know that $g\in L^3(\R)$, $N_g$ is also continuous and we get that $N_g=N_f$ everywhere.
In particular,
$$
\norm{g}_3=N_g(0,0)=N_f(0,0)=\norm{f}_3.
$$
Finally, we have equality in the Cauchy-Schwarz inequality~:
$$
\int g^2=\int g^{1/2}g^{3/2}\leq
\left(\int g\right)^{1/2}\left(\int g^3\right)^{1/2}
=\norm{g}_1^{1/2}\norm{g}_3^{3/2}=\norm{f}_1^{1/2}\norm{f}_3^{3/2}=\norm{f}_2^2
=\norm{g}_2^2.
$$
It follows that $g^{3/2}=\lambda g^{1/2}$ which implies that
$g=\lambda\chi_{\mbox{supp}g}$. Finally, from $\int g^2=\int g$ we get $\lambda=1$
and $g=\chi_{\mbox{supp}g}$.
\end{proof}

\subsection{A restricted problem}
\label{sec:restricted}
One might think that if $f$ and $g$ are characteristic functions
of sets of finite measure it would be impossible to arrange
for examples similar to those of the previous section.

However the situation is subtle, and one should remark that
it is known (\cite{HJ} p.\ 376, \cite{Ka,KV})
that there are measurable sets of finite measure on $\R$
of which the Fourier Transform of their characteristic function vanishes in
a prescribed interval. These sets are of the form
$$
E=\bigcup_{k\in\Z}[k-\lambda_k,k+\lambda_k].
$$
with $0\leq\lambda_k\leq1/2$ and $(\lambda_k)\in\ell^1$.
We will consider slightly more general sets
namely open sets $E = \bigcup_{k\in\Z}(\alpha_k,\beta_k)$, with
$\alpha_k < \beta_k < \alpha_{k+1}$,
such that their complement can be written as the union of closed intervals
whose length is bounded below.
We call such sets {\em open sets with lower bounded gaps}.

We will now show that such sets are characterized up to translations by there $3$-deck.
For this, let us first note that such sets admit the following characterization:

\begin{lemma}
\label{lem:caroswmbg}
Let $E\subset\R$. The following are equivalent:
\begin{enumerate}
\item there exists an $\eps>0$ such that, for every $x,y\in E$ with $|x-y|<\eps$,
the interval $[x,y]$ is contained in $E$;

\item up to a set of measure $0$, $E$ is an open set with lower bounded gaps.
\end{enumerate}
\end{lemma}

\begin{proof} That (2) implies (1) is obvious. Let us assume that (1) holds
and let $Z$ be set of points in $E$ that cannot be approximated
either from the left or from the right by points of $E$.
The set $Z$ is countable as we can assign disjoint open intervals to its points.

Now, let $x\in E\setminus Z$.
By construction of $Z$,
there exists $y\in E$ such that $x<y<x+\eps$. Property (1) then implies that
$[x,y]\subset E$. Similarly there is $z\in E$ with $x-\eps<z<x$ and $[z,x]\subset E$.
It follows that $x\in ]z,y[\subset E$ and $E\setminus Z$ is open.

It is then obvious that $E$ has gaps bounded bellow by $\eps$.
\end{proof}

We are now in position to prove the following:

\begin{theorem}
Let $E$ be an open set of finite measure with lower bounded gaps.
Assume that the set $F$ has the same $3$-deck as $E$.
Then $F$ is a translate of $E$.
\end{theorem}

\begin{proof}
Let us first prove that, up to a set of measure $0$, $F$ is also an open set
with lower bounded gaps.
For this, observe that Property (1) of Lemma \ref{lem:caroswmbg}
holds if and only if there is $\epsilon>0$ such that
$$
G_E(x,y):=\int_\R \chi_E(t)\chi_{{}^cE}(t+x)\chi_E(t+y)dt=0,
\ \ \mbox{whenever $0<x<y<\epsilon$}.
$$
Further $G_E(x,y)=N_E(0,y)-N_E(x,y)$ and, as $N_F=N_E$,
it follows that $F$ is
also an open set with lower bounded gaps (up to measure $0$).

As above, write
$\dst E=\bigcup_{k\in\Z}(\alpha_k,\beta_k)$ with
$\alpha_k<\beta_k<\alpha_{k+1}$ and
$\Gamma_E:=\inf(\alpha_{k+1}-\beta_k)>0$ be the size of the
smallest gap of $E$.
Similarly, let $\Gamma_F$ be the size of the smallest gap of $F$.

Note that $N_E(x,y)=f*g(-x)$ where
$f(t)=\chi_E(t)\chi_E(t+y)$ and $g(t)=\chi_E(-t)$.
It follows that $\partial_xN_E(x,y)=-f*g'(-x)$. Further,
$$
g'=\sum_{k\in\Z}\delta_{-\beta_k}-\delta_{-\alpha_k},
$$
where $\delta_x$ denotes a unit point mass at $x\in\R$,
so that
\begin{align}
\partial_xN_E(x,y)=&-\sum_{k\in\Z}f*\delta_{-\beta_k}(-x)-f*\delta_{-\alpha_k}(-x)\notag\\
&=\sum_{k\in\Z}f(\alpha_k-x)-f(\beta_k-x)\notag\\
&=\sum_{k\in\Z}\chi_E(\alpha_k-x)\chi_E(\alpha_k-x+y)-\chi_E(\beta_k-x)\chi_E(\beta_k-x+y)
\label{parxne}
\end{align}
(the convergence of this series will soon be obvious in the case of interest to us).
Now, if $-\Gamma_E<x<0$, then $\beta_k<\beta_k-x<\alpha_{k+1}$ so that
$\chi_E(\beta_k-x)=0$.
Further, as $\alpha_k<\alpha_k-x<\beta_k+\Gamma_E\leq\alpha_{k+1}$,
$\chi_E(\alpha_k-x)=1$ only when $\beta_k-\alpha_k>x$. As $E$ is of finite measure,
$\beta_k-\alpha_k\to0$ so that the sum (\ref{parxne}) is finite and
\begin{equation}
\label{parxne2}
\partial_xN_E(x,y)=\sum_{k\in\Z/\ \beta_k-\alpha_k>x}\chi_E(\alpha_k-x+y).
\end{equation}
Fix $x$ small enough to have $x<\min(\Gamma_E,\Gamma_F)$ and
let $D_E$ (resp.\ $D_F$) be the Fourier transform in the $y$ variable of
$\partial_xN_E(x,y)$ (resp.\ of $\partial_xN_F(x,y)$).
Then $D_E$ is of the form $D_E(\xi)=\widehat{\chi_E}(\xi)\ffi(\xi)$ where
$\ffi(\xi)=\sum_{k\in\Z/\ \beta_k-\alpha_k>x}e^{2i\pi(\alpha_k-x)\xi}$
is an analytic function, that is not identically $0$ if $x$ is small enough and such that
$\ffi(-\xi)=\overline{\ffi(\xi)}$.
Further, $D_F$ is of the same form $D_F(\xi)=\widehat{\chi_F}(\xi)\psi(\xi)$.

As we assumed $N_F=N_E$, we also have $\partial_x N_E =\partial_x N_F$ and
consequently $D_F=D_E$ so that
$$
D_E(\xi) D_E(\eta) D_E(-\xi-\eta) =
 D_F(\xi) D_F(\eta) D_F(-\xi-\eta).
$$
Using the particular form of $D_E$ and $D_F$, this is
$$
\widehat{\chi_E}(\xi)\widehat{\chi_E}(\eta)\widehat{\chi_E}(-\xi-\eta)
 \ffi(\xi)\ffi(\eta)\ffi(-\xi-\eta) =
\widehat{\chi_E}(\xi)\widehat{\chi_E}(\eta)\widehat{\chi_E}(-\xi-\eta)
 \psi(\xi)\psi(\eta)\psi(-\xi-\eta),
$$
for all $\xi, \eta \in \R$.
But
$$
\widehat{\chi_E}(\xi)\widehat{\chi_E}(\eta)
\widehat{\chi_E}(-\xi-\eta)=
\widehat{\chi_F}(\xi)\widehat{\chi_F}(\eta)
\widehat{\chi_F}(-\xi-\eta)
$$
and is not zero if $\xi,\eta$ are in some neighborhood of $0$ so that
\begin{equation}
\ffi(\xi)\ffi(\eta)\ffi(-\xi-\eta)=\psi(\xi)\psi(\eta)\psi(-\xi-\eta)
\label{eqffipsi}
\end{equation}
in that neighborhood. By analyticity of $\ffi,\psi$, (\ref{eqffipsi})
is valid for every $\xi,\eta\in\R$. It is then easy to see that the
function $\gamma=\psi/\ffi$ is of modulus 1 on the real line and
satisfies $\gamma(\xi+\eta)=\gamma(\xi)\gamma(\eta)$ (excepted for $\xi,\eta$ in some discrete set)
so that $\gamma(\xi)=e^{2i\pi a\xi}$ for some $a\in\R$
and $\psi(\xi)=e^{2i\pi a\xi}\ffi(\xi)$ ({\it see} Remark 3 after the proof of Lemma \ref{lem:rk3}).

Returning to $D_E=D_F$, we get $\widehat{\chi_E}(\xi)\ffi(\xi)=\widehat{\chi_F}(\xi)\psi(\xi)=
e^{2i\pi a\xi}\widehat{\chi_F}(\xi)\ffi(\xi)$. As $\ffi(\xi)$ is analytic
and not identically $0$ it follows that $\widehat{\chi_E}(\xi)=e^{2i\pi a\xi}\widehat{\chi_F}(\xi)$
almost everywhere and $F$ is a translate of $E$.
\end{proof}

\section{The $k$-deck problem on $\Z/n\Z$}
\label{sec:cyclic}

\subsection{Introduction}
The aim of this section is essentially to bring together the knowledge
about the problem in the combinatorics and crystallography communities.
After summarizing known results and open questions, we will proceed with
our own contributions.

In order to do so, let us start with some notations: $k,n$ will be two
positive integers, $k\geq2$ and $\zeta=e^{2\pi i/n}$. 
We will not distinguish between sequences of $n$ elements, $n$-periodic
sequences and functions on $\Z/n\Z$. If 
$f=(f_0,f_1,\ldots,f_{n-1})$ is such
a function, its Fourier transform is defined by
$$
\widehat{f}(l)=\dst\sum_{j=0}^{n-1}f_j\zeta^{jl}\quad l\in\Z
$$
with the natural extensions to higher dimension.
The $k$-deck of $f$ is the function on $(\Z/n\Z)^{k-1}$ given by
$$
N_f^k(j_1,\ldots,j_{k-1})=\sum_{j=0}^{k-1}f_jf_{j+j_1}\ldots f_{j+j_{k-1}}.
$$
Its Fourier transform is then given by
$$
\widehat{N_f^k}(l_1,\ldots,l_{k-1})=\widehat{f}(l_1)\ldots\widehat{f}(l_{k-1})
\widehat{f}(-l_1-\ldots-l_{k-1}).
$$
Now let $f,g$ be two functions on $\Z/n\Z$ that take only non-negative values.
As in the previous section, $f,g$ have same $k$-deck
if and only if we can write $\widehat{g}(k)=\xi(k)\widehat{f}(k)$
where $\xi(k)$ is unimodular and such that
\begin{equation}
\label{fundznz}
\xi(l_1+\ldots+l_{k-1})=\xi(l_1)\ldots\xi(l_{k-1})
\end{equation}
whenever $l_1,\ldots,l_{k-1}\in\mbox{supp}\,\widehat{f}$ are such
that $l_1+\ldots+l_{k-1}\in\mbox{supp}\,\widehat{f}$. Note that $\xi$ is 
a priori only defined on $\mbox{supp}\,\widehat{f}$. Our aim here is to 
know whether or not such a $\xi$ extends to a character of $\Z/n\Z$
in which case $g$ is a translate of $f$. This is of course immediate if
$\widehat{f}$ does not cancel and $k\geq3$. An other trivial case
is when $\widehat{f}$ cancels everywhere excepted at $0$, that is
$f$ is constant.

\subsection{A quick overview of the $2$-deck problem}
Let us first concentrate on the case $k=2$ when condition (\ref{fundznz})
is void. We are seeking $g$ such that $|\widehat{g}|=|\widehat{f}|$, which is
a well known as the {\it phase retrieval problem} in crystallography. The
study of this problem originates in the work of Patterson and a full solution
has been given by Rosenblatt (see \cite{Ro},\ \cite{RoS} and the references therein).

Let us briefly describe what may happen in the case $f=\chi_E$ and $g=\chi_F$
with $E,F$ subsets of $\Z/n\Z$. One easily sees that $E$ and $F$ have same
2-deck if and only if $E-E=F-F$ (counted with multiplicity). From this, one
easily gets the trivial solutions $F=E-a$ and $F=-E-a$. If this are the
only solutions, we say that $E$ is uniquely determined up to translations and
reflections by its $2$-deck. 

Note that not all sets are uniquely determined by there $2$-deck.
For instance, one may arrange for two sets $A,B$ of cardinality at least 2
to be such that $A-B$ and $A+B$ (counted with multiplicity) are still sets
({\it i.e.} every element has multiplicity $1$).
One then immediately gets that $E=A-B$ and $F=A+B$ have same $2$-deck, but
as $A$ and $B$ have more than $2$ elements, $F$ is not of the form $E-a$ nor $-E-a$.
Rosenblatt showed that this is almost always the case and gave a general solution to the problem, which is a bit to long to be summarized here. However, as he noticed himself, 
his solution may be difficult to use in practice. For instance, the following questions are open:
\vskip3pt
{\bf Question 1 (\cite{RS1} conjecture 2):} {\sl Does the proportion of subsets
of $\Z/n\Z$ that are not determined up to translation and reflection by their 2-deck go to
$0$ as $n$ goes to $\infty$.}
\vskip3pt
{\bf Question 2:} {\sl How many solutions (up to translations and reflection)
to the $2$-deck problem can a subset of $\Z/n\Z$ of cardinality $k$ have.}
\vskip3pt
This question is implicit in \cite{Ro}.

\subsection{Zeros in the spectrum of an indicator function}

As noted in the introduction of this section, the $k$-deck problem for $k\geq3$ is trivially
solved if $\widehat{\chi_E}$ has no zeros. Before pursuing with our study we
will therefore gather here some information about possible zeros of $\widehat{\chi_E}$.

First assume that $n=pq$ ($p,q$ not necessarily primes) and let $f$ be
$n$-periodic. Assume that $f$ is also $p$-periodic and write
$a=(f_0,\ldots,f_{p-1})$ for the period. We may see $a$ as a function on $\Z/p\Z$
and still write $\widehat{a}$ for its Fourier transform on $\Z/p\Z$.
An immediate computation shows that
$\widehat{f}$ is supported on the subgroup $\{0,q,\ldots,(p-1)q\}$
of $\Z/n\Z$ and that $\widehat{f}(lq)=q\widehat{a}(l)$. In particular, we get:
\vskip3pt
\noindent{\bf Fact 1.} {\sl A subset $E$ of $\Z/p\Z$ is uniquely
determined up to translation by its $k$-deck, if and only if its $q$-periodic extension 
to $\Z/pq\Z$ is uniquely determined up to translation by its $k$-deck.}
\vskip3pt
\begin{notation}
We will write $q\Z/p\Z=\{0,q,\ldots,(p-1)q\}$.
\end{notation}

Recall that an $n$-th root of unity $\zeta^p$, $0\leq p<n$ is said to be
\emph{primitive} if it is not an $m$-th root of unity for some $m<n$,
that is, $p$ and $n$ are relatively prime. In particular, $\zeta^p$ is a
primitive $\frac{n}{(n,p)}$-th root of unity.
The \emph{cyclotomic polynomial} of order $n$ is then defined as
$$
\Phi_n(x)=\prod_{\begin{matrix}\omega\mbox{ primitive}\\ n\mbox{-th roots}\\
\end{matrix}}(x-\omega)=\prod_{\{j\,:\ (j,n)=1\}}(x-\zeta^j).
$$
It is well known that $\Phi_n$ is the minimal polynomial of any
primitive $n$-th root so that if $P$ is a polynomial such that $P(\zeta^l)=0$
then one has a factorization of $P$: $P(x)=\Phi_{n/(n,l)}(x)Q(x)$,
and $P(\zeta^j)=0$ whenever $j$ and $n/(n,l)$ are relatively prime.

For instance, for $E\subset\Z/n\Z$, let $\dst P_E(x)=\sum_{j\in E}x^j$,
so that $P_E(\zeta^l)=\widehat{\chi_E}(l)$. This leads us to:
\vskip3pt
\noindent{\bf Fact 2.} {\sl If $\widehat{\chi_E}(l)=0$ for some $l\not=0$, then
$\widehat{\chi_E}(j)=0$ for all $1\leq j<n$ with $j$ and $n/(n,l)$
relatively prime.}
\vskip3pt
Two particular cases are of interest to us:
\begin{enumerate}
\item if $n=p^a$ with $p$ prime, then only the three following cases may occur:
\begin{enumerate}
\item $\widehat{\chi_E}(l)\not=0$ for all $l$,
\item $\widehat{\chi_E}(l)=0$ only for $l$ of the form $l=qp^{a-b}$, $1\leq q<p^b$,
\item $\widehat{\chi_E}(l)\not=0$ only for $l$ of the form $l=qp^{a-b}$, $1\leq q<p^b$,
in which case $E$ is $p^b$-periodic.
\item $\widehat{\chi_E}(l)=0$ for all $l\not=0$, in which case $E=\Z/n\Z$.
\end{enumerate}
\item If $n=pq$ with $p,q$ two distinct primes, then only the five following 
cases may occur:
\begin{enumerate}
\item $\widehat{\chi_E}(l)\not=0$ for all $l$,
\item $\widehat{\chi_E}(l)\not=0$ unless $l$ is either of the form $l=jq$, $1\leq j<p$, or of the form $k=jp$, $1\leq j<q$, in which case
$\widehat{\chi_E}(l)=\widehat{f_p}+\widehat{f_q}$ with $f_p$ $p$-periodic and $f_q$ 
$q$-periodic.
\item $\widehat{\chi_E}(l)=0$ unless $l=jq$, $1\leq j<p$, in which case
$E$ is $q$-periodic.  
\item $\widehat{\chi_E}(l)=0$ unless $l=jp$, $1\leq j<q$, in which case
$E$ is $p$-periodic.  
\item $\widehat{\chi_E}(l)=0$ for all $l\not=0$, in which case $E=\Z/n\Z$.
\end{enumerate}
\end{enumerate}

\subsection{Known results on $\Z/n\Z$}

Let us now turn to the $k$-deck problem. Most questions that arise naturally 
have been solved either in \cite{RS1} or in \cite{GM}. More precisely, let us summarize what we consider as the main facts:
\begin{enumerate}
\item If $n$ is a prime number, then every set (and even every integer valued
function) is uniquely determined up to translations by its $3$-deck;
\cite{RS1} Theorem 3. The proof in \cite{RS1} is
rather involved and we will give a simpler one bellow by a more explicit use of Fourier analysis.

\item Assume $n$ is odd and $E\subset\Z/n\Z$. Let $a=\chi_E$ then if
$\widehat{a}(1)\not=0$, $a$ is determined up to translations by its $3$-deck;
\cite{GM}, Theorem 2.

\item Every subset of $\Z/n\Z$ is uniquely determined up to translation by its
$4$-deck; \cite{GM} Theorem 5. This is further extended to integer valued
functions in \cite{GM} Theorem 4: they are still determined up to translation by the $4$-deck if $n$ is odd but the $6$-deck is needed when $n$ is even.
Moreover, examples are given there to show that these results are optimal.
This disproves Conjecture 1 of \cite{RS1}. These examples show that as soon as $n$ has $3$ factors, $n=pqr$
with $p\not=q$ prime, then
there are sets in $\Z/n\Z$ that are not uniquely determined
up to translations by there $3$-deck.

\item The proportion of subsets of $\Z/n\Z$ that are not uniquely determined
up to translations by there $3$-deck goes to $0$ as $n\to\infty$;
\cite{RS1} Theorem 4. This is done by proving that the proportion of subsets 
$E$ of $\Z/n\Z$ such that $\widehat{\chi_E}$ has a zero is at most 
$O(n^{-1/2+\eps})$.
\end{enumerate}

\subsection{Some new results}

We will now settle the two remaining cases in which the $3$-deck may suffice: $n=p^a$ with $p$ prime
and $n=pq$ with $p\not=q$ prime.

\begin{theorem}
Let $p$ be a prime number, $p\geq3$ and let $n=p^a$. Then every subset $E$ off $\Z/n\Z$ is uniquely determined up to translations by its $3$-deck.
\end{theorem}

\begin{proof}
Let us first prove the case $n=p$ a prime. In this case, the particular cases of
Fact 2 imply that either:
\begin{itemize}
\item[---] for all $l$, $\widehat{\chi_E}(l)\not=0$, or
\item[---] $E=\Z/n\Z$.
\end{itemize}
In both cases, $E$ is uniquely determined up to translations by its $3$-deck.
This simplifies the proof in \cite{RS1}.

Assume now that $p\geq3$ and that we have proved that every subset of $\Z/p^a\Z$ is 
uniquely determined up to translations by its $3$-deck and let $E\subset\Z/p^{a+1}\Z$.
Again, there are only four cases that may happen:
\begin{itemize}
\item[---] for all $l$, $\widehat{\chi_E}(l)\not=0$, or
\item[---] $E=\Z/n\Z$, or
\item[---] $E$ is $p^b$-periodic for some $1\leq b\leq a$, or
\item[---] $\widehat{\chi_E}(l)=0$ only for $l=qp^b$, $1\leq q<p^{a-b}$.
\end{itemize}
The two first cases are trivial whereas the third case is settled using the
induction hypothesis and Fact 1. Let us now assume we are in the last case.

Assume that $g$ has same $3$-deck as $\chi_E$ and write 
$\widehat{g}(l)=\xi(l)\widehat{\chi_E}(l)$ where $\xi$ satisfies
\begin{equation}
\xi(l_1+l_2)=\xi(l_1)\xi(l_2)\qquad\mbox{for all }l_1,l_2\in\mbox{supp}\,\widehat{\chi_E}
\mbox{ such that }l_1+l_2\in\mbox{supp}\,\widehat{\chi_E}.
\label{endznz}
\end{equation}
As $p\geq3$, our assumption on $E$ implies $\widehat{\chi_E}(1)\not=0$ and  
$\widehat{\chi_E}(2)\not=0$, so that (\ref{endznz}) implies
$\xi(2)=\xi(1+1)=\xi(1)^2$.

Next, assume that we have proved that $\xi(l)=\xi(1)^l$ for all 
$l\in\mbox{supp}\,\widehat{\chi_E}\cap[0,l_0]$. We are then in one of the following cases:
\begin{itemize}
\item[---] $l_0,l_0+1\in\mbox{supp}\,\widehat{\chi_E}$ and then (\ref{endznz}) implies
$\xi(l_0+1)=\xi(l_0)\xi(1)=\xi(1)^{l_0+1}$ by the induction hypothesis.
\item[---] $l_0+1\notin\mbox{supp}\,\widehat{\chi_E}$ and there is nothing to prove.
\item[---] $l_0\notin\mbox{supp}\,\widehat{\chi_E}$ but $l_0+1\in\mbox{supp}\,\widehat{\chi_E}$. But then $l_0=qp^b$ for some $1\leq q<p^{a-b}$ and as $p\geq3$,
$l_0-1$ is not of that form and so $l_0-1\in\mbox{supp}\,\widehat{\chi_E}$.
But then (\ref{endznz}) implies
$\xi(l_0+1)=\xi(l_0-1)\xi(2)=\xi(1)^{l_0+1}$.
\end{itemize}
We conclude that $\xi$ is a character of $\Z/n\Z$ and that $g$ is a translate of $\chi_E$.
\end{proof}

\begin{remark} The case $p=2$ can be settled similarly by using the $4$-deck instead of 
the $3$-deck. One may then write $\xi(k_0+1)=\xi(k_0-1)\xi(1)\xi(1)$,
but this is not stronger then the result in \cite{GM}.

Note that the proof also holds for rational-valued functions on $\Z/n\Z$.
\end{remark}

\begin{theorem} Let $n=pq$ with $p,q$ two distinct primes. Then every subset of $\Z/n\Z$ is uniquely determined up to translations by its $3$-deck.
\end{theorem}

\begin{proof}
In the case $n=pq$, four cases may happen:
\begin{itemize}
\item[---] for all $l$, $\widehat{\chi_E}(l)\not=0$, or
\item[---] $E=\Z/n\Z$, or
\item[---] $E$ is $p$ or $q$ -periodic, or
\item[---] $\widehat{\chi_E}(l)=0$ if and only if $l$ is neither a multiple of $p$ nor a multiple of $q$.
\end{itemize}
The two first cases are again trivial whereas the third one is again settled using Fact 1
and the case $n$ prime. 

Let us now exclude the last case using the fact that $\widehat{\chi_E}=\widehat{\chi_E}*\widehat{\chi_E}$.

Define $f_p$, $f_q$ by
$$
\widehat{f_p}(k)=\begin{cases}
\frac{1}{2}\widehat{\chi_E}(0)&\mbox{if }l=0\\
\widehat{\chi_E}(l)&\mbox{$l\not=0$ a multiple of $p$}\\
0&\mbox{else}\\
\end{cases},\quad
\widehat{f_q}(l)=\begin{cases}
\frac{1}{2}\widehat{\chi_E}(0)&\mbox{if }l=0\\
\widehat{\chi_E}(l)&\mbox{$l\not=0$ a multiple of $q$}\\
0&\mbox{else}\\
\end{cases},
$$
so that $\chi_E=f_p+f_q$. By assumption, $\widehat{f_p}$ (resp. $\widehat{f_q}$) has support $\{0,q,\ldots,(p-1)q\}$ (resp. $\{0,p,\ldots,(q-1)p\}$). But $\widehat{\chi_E}=\widehat{\chi_E}*\widehat{\chi_E}$, that is 
$\widehat{f_p}*\widehat{f_p}+\widehat{f_q}*\widehat{f_q}+2\widehat{f_p}*\widehat{f_q}=
\widehat{f_p}+\widehat{f_q}$. As $\widehat{f_p}*\widehat{f_p}$
is supported in $q\Z/p\Z+q\Z/p\Z=q\Z/p\Z$,
$\widehat{f_q}*\widehat{f_q}$ is supported in $p\Z/q\Z$, we get that
$\widehat{f_p}*\widehat{f_q}$ has to be supported in $q\Z/p\Z\cup p\Z/q\Z$.
Now let $j\in\Z/n\Z$ be outside this set. There exists a unique $l_0\in\{0,\ldots,p-1\}$ and a unique $j_0\in\{0,\ldots,q-1\}$
such that $j=l_0q+j_0p$. Further
$$
0=\widehat{f_p}*\widehat{f_q}(j)=\sum_{l=0}^{p-1}\widehat{f_p}(lq)\widehat{f_q}(j-lq)
=\widehat{f_p}(l_0q)\widehat{f_q}(l_0p)
$$
so that one of $\widehat{f_p}(l_0q)$ $\widehat{f_q}(l_0p)$ is zero. This contradicts the assumption on the support of these functions.
\end{proof}

\begin{remark} In the case of a rational valued function $f$, again the same four cases
are to be considered and only the last one causes problems. If $g$ has same $3$-deck as 
$f$, we may write $f=f_p+f_q$ and $g=g_p+g_q$ with $f_p$ , $f_q$, $g_p$, $g_q$ defined as above replacing $\chi_E$ by $f$ or $g$.

It follows that $f_p$, $g_p$ (resp. $f_q$, $g_q$) are $p$-periodic (resp. $q$-periodic)
functions, their Fourier transforms over $\Z/p\Z$ (resp. $\Z/q\Z$) don't vanish
and they have same $3$-deck. Therefore $g_p$, (resp. $g_q$) is a translate of
$f_p$, (resp. $f_q$) on $\Z/p\Z$ (resp. $\Z/q\Z$):
$$
\widehat{g_p}(l)=e^{2i\pi lj/pq}\widehat{f_p}(l)
\ \mbox{and}\ 
\widehat{g_q}(l)=e^{2i\pi ll/pq}\widehat{f_q}(l)
$$
so that $g(x)=f_p(x-j)+f_q(x-l)$ and this are the only possible solutions of the $3$-deck problem for $f$.
\end{remark}

Let us finally note that Grunbaum and Moore proved that one can not do any better. Indeed, if $n=pqr$ with $p,q$ two distinct primes and $r\geq3$, let $A=qr\Z/p\Z$ and $B=pr\Z/q\Z$. Let $E=A\cup(B+1)$ and $F=A\cup(B+2)$, so that
$E$ and $F$ are not translates of each other. Further, 
$\chi_E(x)=\chi_A(x)+\chi_B(x-1)$, $\chi_F=\chi_A(x)+\chi_B(x-2)$ so that 
$\widehat{\chi_E}(l)=\widehat{\chi_A}(l)+e^{2il\pi/n}\widehat{\chi_B}(l)$ and 
$\widehat{\chi_F}(l)=\widehat{\chi_A}(l)+e^{4il\pi/n}\widehat{\chi_B}(l)$.

Finally, as $\widehat{\chi_A}$ is supported in $qr\Z/p\Z$
while $\widehat{\chi_B}$ is supported in $pr\Z/q\Z$, it is easy to
see that $E$ and $F$ have same $3$-deck.

Note that if $r$ is also prime, this are essentially the only sets that have same $3$-deck
but are not translates of each other.


\appendix

\section{Continuity properties of $N_f$}
Let $k$ be an integer and for $f_0,\ldots,f_k$ be simple functions
on $\R$, set $x_0=0$ and define
$$
N_{f_0,\ldots,f_k}(x_1,\ldots,x_k)=
\int\prod_{j=0}^k f_j(t-x_j)dt.
$$
If $f_0=\cdots=f_k$, we simply write $N_f$ for $N_{f,\ldots,f}$.

\begin{lemma}
\label{lem:boundkcor}
Let $p_1,\ldots,p_k\geq1$. Then
$$
\norm{N_{f_0,\ldots,f_k}}_{L^r(\R^k)}\leq \prod_{j=0}^k\norm{f_j}_{L^{p_j}(\R)}
$$
holds for every measurable step functions on $\R$ if and only
if $\dst 1+\frac{k}{r}=\sum_{j=0}^k\frac{1}{p_j}$.

If this last condition is satisfied, we may define $N_{f_0,\ldots,f_k}$
for $f_j\in L^{p_j}$ by a usual approximation process.
\end{lemma}

\begin{proof} Replacing $f_j(t)$ by $f_j(\lambda t)$ and letting $\lambda$ go
to $0$ and to $+\infty$ shows that the condition
$\dst 1+\frac{k}{r}=\sum_{j=0}^k\frac{1}{p_j}$ is necessary.

For the converse, assume that $\dst 1+\frac{k}{r}=\sum_{j=0}^k\frac{1}{p_j}$
and note that as not all $p_j$'s are $1$, $r>1$.
Let $\ffi_1,\ldots\ffi_k\in L^{r'}(\R)$ with $\dst\frac{1}{r}+\frac{1}{r'}=1$.
By density of functions of the form $\ffi_1(x_1)\ldots\ffi_k(x_k)$ in $L^{r'}(\R)$ and by duality,
it is enough to prove that
$$
\int_{\R^k}\ffi_1(x_1)\cdots\ffi_k(x_k)N_{f_0,\ldots,f_k}(x_1,\ldots,x_k)dx_1\cdots,dx_k
\leq \norm{\ffi_i}_{s'}\prod\norm{f_j}_{p_j}.
$$
But, by Fubini
\begin{align}
\int_{\R^k}\ffi_1(x_1)\cdots\ffi_k(x_k)&N_{f_0,\ldots,f_k}(x_1,\ldots,x_k)dx_1\cdots,dx_k\notag\\
=&\int_\R\int_{\R^k}\ffi_1(x_1)\cdots\ffi_k(x_k)f_0(t)f_1(t-x_1)\cdots f_k(t-x_k)dtdx_1\cdots dx_k\notag\\
=&\int_\R f_0(t)\int_\R\ffi_1(x_1)f_1(t-x_1)dx_1\cdots\ffi_k(x_k)f_1(t-x_k)dx_k\notag\\
=&\int_\R f_0(t) \ffi_1*\tilde f_1(t)\cdots\ffi_k*\tilde f_k(t)dt,\label{eq:estkdeck}
\end{align}
with 
Now, define $\tilde p_j$ by $\dst1+\frac{1}{\tilde p_j}=\frac{1}{p_j}+\frac{1}{r'}=\frac{1}{p_j}+1-\frac{1}{r}$
and note that
$$
\dst\frac{1}{p_0}+\frac{1}{\tilde p_1}+\cdots+\frac{1}{\tilde p_k}=
\frac{1}{p_0}+\frac{1}{p_1}+\cdots+\frac{1}{p_k}-\frac{k}{r}=1,
$$
so that H\"older's inequality implies that (\ref{eq:estkdeck}) is
$$
\leq \norm{f_0}_{p_0}\norm{\ffi_1*\tilde f_1}_{\tilde p_1}\cdots\norm{\ffi_k*\tilde f_k}_{\tilde p_k}.
$$
To conclude, as $1+\frac{1}{\tilde p_j}=\frac{1}{p_j}+\frac{1}{r'}$, Young's inequality implies that this is
$$
\leq 
\norm{f_0}_{p_0}\norm{\ffi_1}_{r'}\norm{f_1}{p_1}\cdots\norm{\ffi_k}_{r'}\norm{f_1}{p_k}
$$
as wanted.
\end{proof}

\begin{lemma} Let $f$ be a non-negative $L^1$-function, then
$$
\lim\limits_{x_1,\ldots,x_k\to 0}N_{f}(x_1,\ldots,x_k)=\int_\R f^k dt.
$$
\end{lemma}

\begin{proof} For sake of simplicity, we will only prove this for $k=2$. So, let
$$
N_f(x,y)=\int_\R f(t)f(t-x)f(t-y)dt.
$$
By H\"older's inequality,
$$
N_f(x,y)\leq\int_\R f^3(t)dt
$$
so that we only have to prove
$\liminf N_f\geq \int_\R f^3(t)dt$.
For this, let us consider the cake representation of $f$:
$$
f(t)=\int_0^\infty\chi_{\{f>u\}}(t)du,
$$
and write (with Fubini)
\begin{align}
N_f(x,y)=&\int_\R\int_0^\infty\int_0^\infty\int_0^\infty
\chi_{\{f>u\}}(t)\chi_{\{f>v\}}(t-x)\chi_{\{f>w\}}(t-y)
dudvdwdt\notag\\
=&\int_0^\infty\int_0^\infty\int_0^\infty
N_{\chi_{\{f>u\}},\chi_{\{f>v\}},\chi_{\{f>w\}}}(x,y)dudvdw.\notag
\end{align}
It follows from Fatou's Lemma that it is enough to prove that
$$
\liminf_{x,y\to0}N_{\chi_E,\chi_F,\chi_G}(x,y)\geq N_{\chi_E,\chi_F,\chi_G}(0,0),
$$
for every measurable sets $E,F,G$ of finite measure.

Now, fix $\delta>0$, by regularity of the Lebesgue measure, there exists open sets
$U,V$ such that $F\subset U$, $\abs{U\setminus F}<\delta$ and
$G\subset V$, $\abs{V\setminus G}<\delta$. Then
\begin{align}
N_{\chi_E,\chi_F,\chi_G}(x,y)=&\int_\R \chi_E(t)\chi_F(t-x)\chi_G(t-y)dt\notag\\
\geq&\int_\R \chi_E(t)\chi_U(t-x)\chi_V(t-y)dt-2\delta\notag\\
\geq\int_\R \chi_E(t)\chi_F(t)\chi_G(t)\chi_U(t-x)\chi_V(t-y)dt-2\delta.\label{eq:contneen0}
\end{align}
But, as $U,V$ are open,
$$
\chi_E(t)\chi_F(t)\chi_G(t)\chi_U(t-x)\chi_V(t-y)\to\chi_E(t)\chi_F(t)\chi_G(t)\chi_U(t)\chi_V(t)
=\chi_E(t)\chi_F(t)\chi_G(t)
$$
when $x,y\to0$ and $\chi_E(t)\chi_F(t)\chi_G(t)\chi_U(t-x)\chi_V(t-y)\leq\chi_E\in L^1(\R)$. By dominated convergence,
we then have
$$
\int_\R \chi_E(t)\chi_F(t)\chi_G(t)\chi_U(t-x)\chi_V(t-y)dt\to\int_\R\chi_E(t)\chi_F(t)\chi_G(t)=N_{\chi_E,\chi_F,\chi_G}(0,0).
$$
As $\delta$ is arbitrary, we get the desired result from (\ref{eq:contneen0}).
\end{proof}

\begin{thebibliography}{99}

\bibitem[AK]{AK}
\textsc{R. L. Adler and A. G. Konheim}
\newblock{\em A note on translation invariants},
Proc. AMS {\bf 13} (1962), 425--428;
{\bf MR 46 \#4172}

\bibitem[BM]{BM}
\textsc{R. H. T. Bates and D. Mnyama}
\newblock {\em The status of practical Fourier phase retrieval},
Advances in Electronics and Electron physics  {\bf 67} (1986), 1--64;
{\bf MR 58 \#372}

\bibitem[Bo]{Bo}
\textsc{J.A. Bondy}
\newblock {\em A graph reconstructor's manual}, in {\it Handbook of
combinatorics}, {\bf Vol. 1, 2}, 3--110, Elsevier,
Amsterdam, 1995; {\bf MR 97a:05129}

\bibitem[BH]{BH}
\textsc{J.A. Bondy and R.L. Hemminger}
\newblock {\em Graph reconstruction---a survey},
J. Graph Theory  {\bf 1} (1977), no.~3, 227--268;
{\bf MR 58 \#372}

\bibitem[CW]{CW}
\textsc{D. Chazan and B. Weiss}
\newblock {\em Higher order autocorrelation functions as translation invariants},
Information and Control {\bf 16} (1970), 378-383;
{\bf MR 43\#1749}

\bibitem[GM]{GM}
\textsc{F.A. Gr\"unbaum and C.C. Moore}
\newblock {\em The use of higher-order invariants in the
determination of generalized Patterson cyclotomic sets},
Acta Cryst. Sect. A {\bf 51} (1995), no.~3, 310--323;
{\bf MR 96e:82087}


\bibitem[HJ]{HJ}
\textsc{V. Havin and B. J\" oricke}
\newblock {\em The uncertainty principle in harmonic analysis},
Ergebnisse der Mathematik und ihrer Grenzgebiete (3), 28;
Springer-Verlag, Berlin, 1994;
{\bf MR 96c:42001}

\bibitem[Hu]{Hu}
\textsc{N. E. Hurt}
\newblock {\em Phase Retrieval and Zero Crossing (Mathematical Methods
in Image Reconstruction)}.
Math. and Its Appl. Kluwer Academic Publisher, 1989.
{\bf MR 92k:94002}

\bibitem[Ja]{Ja}
\textsc{P. Jaming}
\newblock {\em Phase retrieval techniques for radar ambiguity problems},
J. Fourier Anal. Appl.  {\bf 5} (1999), no.~4, 309--329;
{\bf MR 2000g:94007}

\bibitem[Ka]{Ka}
\textsc{P.P. Kargaev}
\newblock {\em The Fourier Transform of the characteristic
function of a set that is vanishing on the interval},
(Russian), Mat. Sb. (N.S.) {\bf 117(159)} (1982), no. 3, 397--411, 432;
{\bf MR 83f:42010}

\bibitem[KV]{KV}
\textsc{P.P. Kargaev and A.L. Volberg}
\newblock {\em Three results concerning the support of functions and their
Fourier Transforms}, Indiana Univ.\ Math.\ J. {\bf 41} (1992), no. 4, 1143--1164;
{\bf MR 94d:42018}

\bibitem[KST]{KST}
\textsc{M. V. Klibanov, P.E. Sacks and A.V. Tikhonravov}
\newblock {\em The phase retrieval problem},
Inverse problems {\bf 11} (1995), no.~1, 1--28;
{\bf MR 95m:35203}

\bibitem[Ko]{Ko}
\textsc{P. Koosis},
\newblock {\em The logarithmic integral, I},
Corrected reprint of the 1988 original;
Cambridge Tracts in Advanced Mathematics, 12;
Cambridge Univ.\ Press, 1998;
{\bf MR 99j:30001}


\bibitem[LL]{LL}
\textsc{Y.T. Lam and K.H. Leung}
\newblock {\em On vanishing sums of roots of unity},
J. Algebra {\bf 224} (2000), no.~1, 91--109;
{\bf MR 2001f:11135}

\bibitem[Mi]{Mi}
\textsc{R. P. Millane}
\newblock {\em Phase retrieval in crystallography and optics},
J. Opt. Soc. Am. A. {\bf 7} (1990), no.~3, 394--411

\bibitem[RS1]{RS1}
\textsc{A.J. Radcliffe and A.D. Scott}
\newblock {\em Reconstructing subsets of $\Z_n$},
J. Combin.\ Theory Ser.\ A {\bf 83} (1998), no.~2, 169--187;
{\bf MR 99k:05158}

\bibitem[RS2]{RS2}
\textsc{A.J. Radcliffe and A.D. Scott}
\newblock {\em Reconstructing subsets of reals},
Electron.\ J. Combin. {\bf 6} (1999), no.~1, Research Paper 20, 7 pp. (electronic);
{\bf MR 2000d:05081}

\bibitem[RT]{RT}
\textsc{D. Rautenbach and E. Triesch}
\newblock {\em A note on the reconstruction of sets of finite measure},
manuscript 2002.

\bibitem[Ro]{Ro}
\textsc{J. Rosenblatt}
\newblock {\em Phase retrieval},
Comm. Math. Phys. {\bf 95} (1984), no.~3, 317--343;
{\bf MR 86k:82075}

\bibitem[RoS]{RoS}
\textsc{J. Rosenblatt and P. Seymour}
\newblock {\em The structure of homometric sets},
SIAM J. Algebraic Discrete Methods {\bf 3} (1982), no.~3, 343--350;
{\bf MR 83m:12029}


\bibitem[Rot]{Rot}
\textsc{J. Rothman}
\newblock {\em Autocorrelation functions as translation invariants in $L^1$ and $L^2$},
J. Fourier Anal. Appl. {\bf 2} (1996), no.~3, 217--225;
{\bf MR 97a43005}

\end{thebibliography}
\end{document}